\documentclass[a4paper,11pt,oneside]{article}
\usepackage{latexsym,amsfonts,amsmath,amssymb,verbatim}

\newcommand{\Z}{\mathbb{Z}}

\newcommand{\C}{\mathbb{C}}
\newcommand{\pres}[2]{\langle\ {#1}\ |\ {#2}\ \rangle}

\newcommand{\tr}{\mathrm{tr}}
\newcommand{\ol}[1]{\bar{#1}}
\newtheorem{theorem}{Theorem}
\newtheorem{lemma}[theorem]{Lemma}

\newenvironment{maintheorem}{\medskip \noindent\textbf{Main Theorem.}
\em}{\medskip}
\newenvironment{proof}{\noindent\textbf{Proof}\newline}{$\square$}
\def\PSL{{PSL}(2,\mathbb{C})}
\def\SL{{SL}(2,\mathbb{C})}
\def\Z{\mathbb{Z}}

\begin{document}
\title{The Tits alternative for generalized triangle groups of type
$(3,4,2)$} \author{James Howie and Gerald Williams} \maketitle
\begin{abstract}
A generalized triangle group is a group that can be presented in
the form
\( G = \pres{x,y}{x^p=y^q=w(x,y)^r=1} \)
where $p,q,r\geq 2$ and $w(x,y)$ is a cyclically reduced word of
length at least $2$ in the free product
$\Z_p*\Z_q=\pres{x,y}{x^p=y^q=1}$.
Rosenberger has conjectured that every generalized triangle group
$G$ satisfies the Tits alternative. It is known that the conjecture
holds except possibly when the triple $(p,q,r)$ is one of $(2,3,2),\
(2,4,2),\ (2,5,2),\ (3,3,2),\ (3,4,2),$ or $(3,5,2)$. In this paper
we show that the Tits alternative holds in the case
$(p,q,r)=(3,4,2)$.
\end{abstract}
\section{Introduction}
A \em generalized triangle group \em is a group that can be
presented in the form
\[ G = \pres{x,y}{x^p=y^q=w(x,y)^r=1} \]
where $p,q,r\geq 2$ and $w(x,y)$  is a cyclically reduced word of
length at least $2$ in the free product
$\Z_p*\Z_q=\pres{x,y}{x^p=y^q=1}$ that is not a proper power. It was
conjectured by Rosenberger~\cite{Rosenberger} that every generalized
triangle group $G$ satisfies the Tits alternative. That is, $G$
either contains a non-abelian free subgroup or has a soluble
subgroup of finite index.

It is now known that the Tits alternative holds for a generalized
triangle group $G$ except possibly when the triple $(p,q,r)$ is one
of $(2,3,2)$, $(2,4,2)$, $(2,5,2)$, $(3,3,2)$, $(3,4,2)$, or
$(3,5,2)$. (See~\cite{BMS} for the case $1/p+1/q+1/r<1$;
see~\cite{FLR} for the case $r\geq 3$; see~\cite{Howie} for the cases $(3,6,2)$ and $(4,4,2)$; and
see~\cite{BKB},\cite{BK},\cite{BK2},\cite{BKB2},\cite{HW} for the
cases $(2,q,2)$ where $q\geq 6$. Many of these results are described
in the survey article~\cite{FRR}.) In this paper we show that the
Tits alternative holds in the case $(3,4,2)$:

\begin{maintheorem}\label{maintheorem}
Let $\Gamma = \pres{x,y}{x^3=y^4=w(x,y)^2=1}$ where $w(x,y) =
x^{\alpha_1}y^{\beta_1} \ldots$ $x^{\alpha_k}y^{\beta_k}$, $1\leq
\alpha_i \leq 2$, $1 \leq \beta_i \leq 3$ for each $1\leq i \leq k$
where $k\geq 1$. Then the Tits alternative holds for $\Gamma$.
\end{maintheorem}

Benyash-Krivets and Barkovich~\cite{BKB2} have proved this result
when $k$ is even, and for this reason we focus on the case when $k$
is odd.

\section{Preliminaries}
We first recall some definitions and well-known facts concerning
generalized triangle groups; further details are available in (for
example)~\cite{FRR}.

Let $G=\pres{x,y}{x^\ell=y^m=w(x,y)^2=1}$ where $w(x,y) =
x^{\alpha_1}y^{\beta_1} \ldots x^{\alpha_k}y^{\beta_k}$, $1\leq
\alpha_i < \ell$, $1 \leq \beta_i < m$ for each $1\leq i \leq k$
where $k\geq 1$. A homomorphism $\rho: G \rightarrow H$ (for some
group $H$) is said to be \em essential \em if
$\rho(x),\rho(y),\rho(w)$ are of orders $\ell,m,2$ respectively.
By~\cite{BMS} $G$ admits an essential representation into $\PSL$.

A projective matrix $A\in \PSL$ is of order $n$ if and only if $\tr
(A) = 2\cos(q\pi/n)$ for some $(q,n)=1$. Note that in $\PSL$ traces
are only defined up to sign. A subgroup of $\PSL$ is said to be \em
elementary \em if it has a soluble subgroup of finite index, and is
said to be \em non-elementary \em otherwise.

Let $\rho: \pres{x,y}{x^\ell=y^m=1} \rightarrow \PSL$ be given by
$x\mapsto X,\ y\mapsto Y$ where $X,Y$ have orders $\ell,m,$
respectively. Then $w(x,y) \mapsto w(X,Y)$. By
Horowitz~\cite{Horowitz} $\tr w(X,Y)$ is a polynomial with integer
coefficients in $\tr X, \tr Y, \tr XY$, of degree $k$ in $\tr XY$.
Since $X,Y$ have orders $\ell,m,$ respectively, we may assume (by
composing $\rho$ with an automorphism of $\pres{x,y}{x^\ell=y^m=1}$,
if necessary), that $\tr X=2\cos(\pi/\ell),\ \tr Y = 2 \cos(\pi/m)$.
Moreover (again by~\cite{Horowitz}) $X$ and $Y$ can be any elements
of $\PSL$ with these traces. We refer to $\tr w(X,Y)$ as the \em
trace polynomial \em of $G$. The representation $\rho$ induces an
essential representation $G \rightarrow \PSL$ if and only if $\tr
\rho (w) = 0$; that is, if and only if $\tr XY$ is a root of $\tr
w(X,Y)$. By~\cite{Horowitz} the leading coefficient of $\tr w(X,Y)$
is given by
\begin{eqnarray} c = \prod_{i=1}^k
\frac{\sin (\alpha_i\pi/\ell)\sin
(\beta_i\pi/m)}{\sin(\pi/\ell)\sin(\pi/m)} .
\label{eq:leadingcoefficient}
\end{eqnarray}
%\begin{eqnarray} c = \prod_{i=1}^k \left(
%\frac{\zeta^{\alpha_i}-\zeta^{-\alpha_i}}{\zeta-\zeta^{-1}} \right)
%\left( \frac{\eta^{\beta_i}-\eta^{-\beta_i}}{\eta-\eta^{-1}}
%\right)\label{eq:leadingcoefficient}
%\end{eqnarray}
%
%where $\zeta = e^{i\pi/\ell},\eta=e^{i\pi/m}$.

Now if $X,Y$ generate a non-elementary subgroup of $\PSL$ then $\rho
(G)$ (and hence $G$) contains a non-abelian free subgroup. Thus in
proving that $G$ contains a non-abelian free subgroup we may assume
that $X,Y$ generate an elementary subgroup of $\PSL$. By
Corollary~2.4 of~\cite{Rosenberger} there are then three
possibilities: (i) $X,Y$ generate a finite subgroup of $\PSL$; (ii)
$\tr [X,Y] = 2$; or (iii) $\tr XY =0$.
The finite subgroups of $\PSL$ are the alternating groups $A_4$ and
$A_5$, the symmetric group $S_4$, cyclic and dihedral groups (see
for example~\cite{CM}). The Fricke identity
\[ \tr [X,Y] = (\tr X)^2 + (\tr Y)^2 + (\tr XY)^2 - (\tr X) (\tr Y) (\tr
XY) -2 \]
implies that (ii) is equivalent to $\tr XY = 2 \cos
(\pi/\ell\pm\pi/m)$. These values occur as roots of $\tr w(X,Y)$ if
and only if $G$ admits an essential cyclic representation. Such a
representation can be realized as $x \mapsto A, y\mapsto B$ where
\[ A = \begin{pmatrix} e^{i\pi/\ell} & 0\\ 0 & e^{-i\pi/\ell} \end{pmatrix}, \quad B =
\begin{pmatrix} e^{\pm i\pi/m} & 0\\ 0 & e^{\mp i\pi/m} \end{pmatrix}. \]
We summarize the above as
\begin{lemma}\label{lem:RosenbergerCor2.4}
Let $G=\pres{x,y}{x^\ell=y^m=w(x,y)^2=1}$. Suppose $G \rightarrow
\PSL$ is an essential representation given by $x \mapsto X, y\mapsto
Y$, where $\tr X =2\cos(\pi/\ell),\ \tr Y = 2 \cos(\pi/m)$. If $G$
does not contain a non-abelian free subgroup then one of the
following occurs:
\begin{enumerate}
\item $X,Y$ generate $A_4,S_4,A_5$ or a finite dihedral group;
\item $\tr XY = 2 \cos (\pi/\ell\pm\pi/m)$;
\item $\tr XY =0$.
\end{enumerate}
Case (2) occurs if and only if $G$ admits an essential cyclic
representation.
\end{lemma}
\section{Proof of Main Theorem}

Throughout this section $\Gamma$ will be the group defined in the Main Theorem.

\begin{lemma}\label{lem:essentialcyclic}
If $\Gamma$ admits an
essential cyclic representation then $\Gamma$ contains a non-abelian
free subgroup.
\end{lemma}

\begin{proof}
Let $\rho : \Gamma\rightarrow \Z_{12}$ be an essential
representation. Then $K = \mathrm{ker} \rho$ has a deficiency zero
presentation with generators
\begin{alignat*}{4}
a_1 &= yxy^{-1}x^{-1},&\quad a_2&=y^2xy^{-2}x^{-1},&\quad a_3&=y^3xy^{-3}x^{-1},\\
a_4 &= xyxy^{-1}x^{-2},&\quad a_5&= xy^2xy^{-2}x^{-2},&\quad a_6 &=
xy^3xy^{-3}x^{-2},
\end{alignat*}
and with relators
\[ W'(a_i,\ldots,a_6,a_1,\ldots,a_{i-1})
W'(y^2a_iy^2,\ldots,y^2a_6y^2,y^2a_1y^2,\ldots,y^2a_{i-1}y^2) \quad
(1 \leq i \leq 6)\]
where $W'$ is a rewrite of $W$.

Let $S= \{\,[a_i,a_j],\ a_i(y^2a_iy^2)\ (1\leq i,j\leq 6)\,\}$, and
let $L,N$ respectively be the normal closures of $S$ and $S \cup
\{a_6\}$ in $K$. Noting that
\begin{alignat*}{3}
y^2a_1y^2 &= a_3a_2^{-1},&\quad y^2a_2y^2 &= a_2^{-1}, &\quad y^2a_3y^2 &= a_1a_2^{-1},\\
y^2a_4y^2 &= a_2a_6a_5^{-1}a_2^{-1},&\quad y^2a_5y^2 &=
a_2a_5^{-1}a_2^{-1}, &\quad y^2a_6y^2 &= a_2a_4a_5^{-1}a_2^{-1},
\end{alignat*}
we have that $K/L \cong \Z^4$ and $K/N\cong \Z^3$, and hence that
$N/N'\neq 0$.

Let $\phi : K \rightarrow K$ be given by $a_i \mapsto y^2a_iy^2$
($1\leq i \leq 6$). It is clear from the presentation of $K$ that
$\phi$ is an automorphism of $K$; furthermore $\phi (N) = N$. In the
abelian group $K/N$, $\phi(a_i) = y^2a_iy^2 = a_i^{-1}$ ($1 \leq i
\leq 6$). That is, $\phi$ induces the antipodal automorphism
$\alpha\mapsto-\alpha$ on $K/N$. By Corollary~3.2 of~\cite{Howie},
$K$ contains a non-abelian free subgroup.
\end{proof}

\medskip

We will write the trace polynomial of $\Gamma$ as $\tau (\lambda) =
\tr w(X,Y)$, where $\tr (X) = 1$, $\tr(Y) = \sqrt{2}$,
$\lambda=\tr(XY)$. By Lemmas~\ref{lem:RosenbergerCor2.4}
and~\ref{lem:essentialcyclic} we may assume that $\tr XY=0$ or $X,Y$
generate $A_4,S_4,$ or $A_5$. But $Y$ has order $4$ so $X,Y$ cannot
generate $A_4$ or $A_5$. If $X,Y$ generate $S_4$ then the product
$XY$ has order $2$ or $4$ so $\tr XY = 0, \pm \sqrt{2}$. Suppose
$\tr XY= -\sqrt{2}$. It follows from the identity
\[ \tr XY + \tr X^{-1}Y = (\tr X) (\tr Y) \]
that $\tr X^{-1}Y = 2\sqrt{2}$. Replacing $X$ by $X^{-1}$ in
Lemma~\ref{lem:RosenbergerCor2.4} shows that $\Gamma$ contains a
non-abelian free subgroup. Thus we may assume that the only roots
$\lambda=\tr XY$ of $\tau$ are $\lambda=0,\sqrt{2}$.
Using~(\ref{eq:leadingcoefficient}) the leading coefficient of
$\tau$ is given by $c=\pm(\sqrt{2})^\kappa$ where $\kappa$ denotes
the number of values of $i$ for which $\beta_i=2$. Hence $\tau
(\lambda)$ takes the form
\begin{eqnarray}
\tau(\lambda) = (\sqrt{2})^\kappa \lambda^s
(\lambda-\sqrt{2})^{k-s}\label{eq:tau}
\end{eqnarray}
where $s\geq 0$. Moreover, Theorem~2 of~\cite{BKB2} implies that the
Main Theorem holds when $k$ is even, so we may assume that $k$ is
odd.

Let
\[    A = \begin{pmatrix} e^{i\pi/3} & 0  \\ 1 & e^{-i\pi/3}  \end{pmatrix},
\quad B = \begin{pmatrix} e^{i\pi/4} & z  \\ 0 & e^{-i\pi/4}
\end{pmatrix}.\]
Then $\tr A =1,\ \tr B = \sqrt{2}$, $\tr AB = z -
(\sqrt{6}-\sqrt{2})/2$. Consider the representation $\rho:
\pres{x,y}{x^3=y^4=1} \rightarrow \PSL$ given by $x \mapsto A,
y\mapsto B$. Then $\tr \rho (x^{\alpha_1}y^{\beta_1}\ldots
x^{\alpha_k}y^{\beta_k}) = \tau (z-(\sqrt{6}-\sqrt{2})/2)$ whose
constant term (by~(\ref{eq:tau})) is
\[ \pm (\sqrt{2})^\kappa ((\sqrt{6}-\sqrt{2})/2)^s
((\sqrt{6}+\sqrt{2})/2)^{k-s} \]
which simplifies to
\[ \pm (\sqrt{2})^\kappa ((\sqrt{6}+\sqrt{2})/2)^{k-2s}. \]
Now the constant term in $\tr (A^{\alpha_1}B^{\beta_1} \ldots
A^{\alpha_k}B^{\beta_k})$ is equal to
\[ 2 \cos \left( \frac{(4\sum_{i=1}^k \alpha_i + 3\sum_{i=1}^k
\beta_i)\pi}{12} \right). \]
Thus $(\sqrt{2})^\kappa ((\sqrt{6}+\sqrt{2})/2))^{k-2s} = 2 \cos
\left( \frac{(4\sum_{i=1}^k \alpha_i + 3\sum_{i=1}^k
\beta_i)\pi}{12} \right)$ and since $k$ is odd, this only happens if
$\kappa =0$ and $k-2s=\pm 1$. It follows that
\begin{eqnarray}\label{eq:15711}
 4\sum_{i=1}^k \alpha_i + 3\sum_{i=1}^k \beta_i = 1,5,7,11 \
\mathrm{mod}\ 12.
\end{eqnarray}

Since $\kappa=0$ there is no value of $i$ for which $\beta_i=2$ and
hence $\Gamma$ maps homomorphically onto the group
\begin{eqnarray}
 \ol{\Gamma} = \pres{x,y}{x^3=y^2=\ol{w}(x,y)^2=1} \label{eq:322}
\end{eqnarray}
where $\ol{w} (x,y) = x^{\alpha_1}y \ldots x^{\alpha_k}y$. If $\ol{w}$ is
a proper power then $\ol{\Gamma}$ contains a non-abelian free
subgroup by~\cite{BMS}. Thus we may assume that $\ol{w}$ is not a
proper power, and so~(\ref{eq:322}) is a presentation of
$\ol{\Gamma}$ as a generalized triangle group.

We will write the trace polynomial of $\ol{\Gamma}$ as $\sigma (\mu)
= \tr \ol{w}(\ol{X},\ol{Y})$, where $\tr (\ol{X}) = 1$, $\tr(\ol{Y})
= 0$, $\mu=\tr(\ol{X}\ol{Y})$. It follows from~(\ref{eq:15711}) that
$\sum_{i=1}^k \alpha_i \neq 0$ mod 3 so $\ol{\Gamma}$ admits no
essential cyclic representation. By
Lemma~\ref{lem:RosenbergerCor2.4} we may assume that $\mu = 0$ or
$\ol{X},\ol{Y}$ generate $A_4,S_4,A_5$ or a finite dihedral group,
in which case $\ol{X}\ol{Y}$ has order $2,3,4,$ or $5$ and hence
$\mu = 0, \pm1, \pm \sqrt{2}, (\pm1 \pm \sqrt{5})/2$. Moreover
$\ol{X}$ is of order 4 in $\SL$ so $\ol{X}^{-1}=-\ol{X}$ and thus
$\tr(\ol{X}^{-1}\ol{Y}) = - \mu$ and $\tr \ol{w}(\ol{X},\ol{Y}) =
(-1)^k \tr\ol{w}(\ol{X}^{-1},\ol{Y})$, so $\sigma_w(\mu) = \pm
\sigma_w(-\mu)$. Thus $\mu$ and $-\mu$ occur as roots of $\sigma$
with equal multiplicity. By~(\ref{eq:leadingcoefficient}) the
leading coefficient of $\sigma$ is $\pm 1$ so
\[ \sigma (\mu) = \pm \mu^{u_1} (\mu^2-1)^{u_2} (\mu^2-2)^{u_3}
(\mu^2 - (3+\sqrt{5})/2)^{u_4} (\mu^2-(3-\sqrt{5})/2)^{u_5}
\]
where $u_1,u_2,u_3,u_4,u_5 \geq 0$ and $u_1+2u_2+2u_3+2u_4+2u_5=k$.
Since $\tr (\ol{X}\ol{Y})$ is a polynomial with integer coefficients
in $\tr\ol{X} = 1 ,\tr \ol{Y} =0$ we have that $u_5=u_4$ so
\begin{eqnarray}\label{eq:sigma}
\sigma (\mu) = \pm \mu^{u_1} (\mu^2-1)^{u_2} (\mu^2-2)^{u_3}
(\mu^4-3\mu^2+1)^{u_4}
\end{eqnarray}
and $u_1+2u_2+2u_3+4u_4 =k$. Let
\[    \tilde{A} = \begin{pmatrix} e^{i\pi/3} & 0  \\ 1 & e^{-i\pi/3}  \end{pmatrix},
\quad \tilde{B} = \begin{pmatrix} i & z  \\ 0 & -i
\end{pmatrix}.\]
Then $\tr \tilde{A} = 1,\ \tr \tilde{B} = 0,\ \tr \tilde{A}\tilde{B}
= z -\sqrt{3}$. Now the constant term in $\sigma (z-\sqrt{3})$ is
$(-\sqrt{3})^{u_1}\cdot 2^{u_2}$. But the constant term in $\tr
(\tilde{A}^{\alpha_1} \tilde{B} \ldots \tilde{A}^{\alpha_k}
\tilde{B})$ is $2 \cos ((2\sum_{i=1}^k \alpha_i + 3k) \pi /3) =
\pm\sqrt{3}$ so $u_1=1, u_2=0$ and thus $k=1 + 2u_3+4u_4$.

\begin{lemma}\label{lem:repeatedS4}
If $\sqrt{2}$ is a repeated root of $\sigma (\mu)$ then $\Gamma$
contains a non-abelian free subgroup.
\end{lemma}

\begin{proof}
Let $q:\Gamma\to\ol{\Gamma}$ denote the canonical epimorphism. By
hypothesis, there is an essential representation
$\rho:\ol{\Gamma}\to PSL (2,\C[\mu]/(\mu-\sqrt{2})^2)$. Indeed, we
can construct $\rho$ explicitly via:
\[\rho(x)= \left(\begin{array}{rr}e^{i\pi/3} &\mu\\ 0
&e^{-i\pi/3}\end{array}\right), \qquad \rho(y)=
\left(\begin{array}{cc} 0 & -1\\
1&0\end{array}\right).\]
Composing this with the canonical epimorphism
\[ \psi:PSL(2,\C[\mu]/(\mu-\sqrt{2})^2)\to PSL(2,\C[\mu]/(\mu-\sqrt{2}))\cong
\PSL\]
gives an essential representation
$\tilde{\rho}=\psi\circ\rho:\ol{\Gamma} \to \PSL$ with image $S_4$,
corresponding to the root $\sqrt{2}$ of the trace polynomial.

Let $\ol{K}$ denote the kernel of $\tilde{\rho}$, $V$ the kernel of
$\psi$, and $K$ the kernel of the composite map $\tilde{\rho}\circ
q:\Gamma\to \PSL$.  Then $V$ is a complex vector space, since its
elements have the form $\pm (I+(\mu-\sqrt{2})A)$ for various
$2\times 2$ matrices $A$, with multiplication
\[ [\pm (I + (\mu-\sqrt{2})A)][\pm (I + (\mu-\sqrt{2})B)]=\pm (I +
(\mu-\sqrt{2})(A+B)).\]

Now $\ol{K}$ is generated by conjugates of $(xy)^4$ and $\rho((xy)^4)
= -I +(\mu-\sqrt{2})M$ where $M = \begin{pmatrix} 2\sqrt{2} &
-2(1+i\sqrt{3})\\ 2(1-i\sqrt{3}) & -2\sqrt{2} \end{pmatrix}$. Since
$M$ is non-zero, $\ol{K}$ (and hence $K$) maps onto the free abelian
group of rank 1. Let $N$ be a normal subgroup of $K$ such that
$K/N\cong \Z$.

Note that $K$ is the fundamental group of a $2$-dimensional
CW-complex $X$ arising from the given presentation of $\Gamma$. This
complex $X$ has $24$ cells of dimension $0$, $48$ cells of dimension
$1$, and $24(\frac14+\frac13+\frac12)=26$ cells of dimension $2$.
Here, $24/4=6$ of the 2-cells (call them $\alpha_1,\ldots,\alpha_6$,
say) arise from the relator $y^{4}$, $24/3=8$
($\alpha_{7},\ldots,\alpha_{14}$, say) arise from the relator $x^3$,
and $24/2=12$ ($\alpha_{15},\ldots,\alpha_{26}$, say) arise from the
relator $w(x,y)^2$. Moreover, $\alpha_1,\dots,\alpha_{6}$ are
attached by maps which are $2$nd powers.
Let $\widehat{X}$ be the regular covering complex of $X$
corresponding to the normal subgroup $N$ of $K$ and let
$\widehat{\alpha}_i$ denote a lift of the 2-cell $\alpha_i$. Then
each of $\widehat{\alpha}_1,\ldots,\widehat{\alpha}_{6}$ is a
$2$-cell attached by a map which is a $2$nd power.

Let $GF_2$ denote the field of 2 elements. Now
$H_2(\widehat{X},GF_2)$ is a subgroup of the 2-chain group
$C_2(\widehat{X},GF_2)$ and since $K/N$ freely permutes the cells of
$\widehat{X}$, $C_2(\widehat{X},GF_2)$ is a free $GF_2(K/N)$-module
on the basis $\widehat{\alpha}_1,\ldots,\widehat{\alpha}_{26}$. Let
$Q$ be the free $GF_2 (K/N)$-submodule of $C_2(\widehat{X},GF_2)$ of
rank 6 generated by $\widehat{\alpha}_1,\ldots,
\widehat{\alpha}_{6}$. Since these 2-cells are attached by maps
which are $2$nd powers, their boundaries in the 1-chain group
$C_1(\widehat{X},GF_2)$ are zero. Thus $Q$ is a subgroup of
$H_2(\widehat{X},GF_2)$. Since the rank of $Q$ is greater than
$\chi(X) = 2$, Theorem~A of~\cite{Howie} implies that $K$, and hence
$\Gamma$, contains a non-abelian free subgroup
\end{proof}

\begin{lemma}\label{lem:repeatedA5}
If $(1 + \sqrt{5})/2$ is a repeated root of $\sigma (\mu)$ then
$\Gamma$ contains a non-abelian free subgroup.
\end{lemma}

\begin{proof}
The proof is similar to that of Lemma~\ref{lem:repeatedS4}. In this
case $\tilde{\rho}$ has image $A_5$, corresponding to the root
$(1+\sqrt{5})/2$. The complex $X$ has 60 0-cells, 120 1-cells, and
$60(\frac{1}{4}+\frac{1}{3}+\frac{1}{2}) =65$ 2-cells (so $\chi (X)
= 5$). Moreover, $60/4 = 15$ of the 2-cells (call them
$\alpha_1,\ldots,\alpha_{15}$, say) are attached by maps which are
2nd powers. As before, the free $GF_2 (K/N)$-submodule, $Q$, of
$C_2(\widehat{X},GF_2)$ of rank 15 generated by
$\widehat{\alpha}_1,\ldots, \widehat{\alpha}_{15}$ is a subgroup of
$H_2(\widehat{X},GF_2)$. Since the rank of $Q$ is greater than $\chi
(X)$, Theorem~A of~\cite{Howie} again implies that $K$ contains a non-abelian free
subgroup.
\end{proof}

\medskip

By Lemmas~\ref{lem:repeatedS4} and~\ref{lem:repeatedA5} we may
assume $u_3,u_4 \leq 1$ so $k\leq 7$. A computer search reveals that
if $k=3$ or $7$ then there is no word $w(x,y)$ such that $\tau
(\lambda)$ is of the form~(\ref{eq:tau}). If $k=5$ then (up to
cyclic permutation, inversion, and automorphisms of $\pres{x}{x^3}$
and $\pres{y}{y^4}$) the only word $w(x,y)$ with $\tau(\lambda)$ of
the form (\ref{eq:tau}) is $w = xyxyx^2y^3x^2yxy^3$. In this case, a
computer search using GAP~\cite{GAP} shows that $\Gamma$ contains a
subgroup of index 4 which maps onto the free group of rank 2. If
$k=1$ then either $\Gamma=\pres{x,y}{x^3=y^4=(xy)^2=1}$ or $\Gamma =
\pres{x,y}{x^3=y^4=(xy^2)^2=1}$. In the first case $\Gamma\cong
S_4$, and in the second $\Gamma$ can be written as an amalgamated
free product
\[ \Gamma =
\pres{x,y^2}{x^3=y^4=(xy^2)^2=1}\underset{\pres{y^2}{y^4}}{*}\pres{y}{y^4}
\]
in which the amalgamated subgroup has index $3$ in the first factor
and index $2$ in the second, and thus $\Gamma$
contains a non-abelian free subgroup. This completes the proof of
the Main Theorem.

\medskip\noindent
{\bf Author addresses:}
\small
\begin{flushleft}
James Howie\\
School of Mathematical and Computer
Sciences\\
Heriot-Watt University\\
Edinburgh EH14 4AS\\
{\tt J.Howie@hw.ac.uk}
\end{flushleft}

\begin{flushleft}
Gerald Williams\\
Institute of Mathematics, Statistics and Actuarial Science\\
University of Kent\\
Canterbury\\
Kent CT2 7NF\\
{\tt g.williams@kent.ac.uk }
\end{flushleft}

\end{document}